\documentclass[leqno,12pt]{amsart}
\usepackage{amssymb} 
\theoremstyle{plain} 
\newtheorem{theorem}{\indent\sc Theorem}[section] 

\newtheorem{proposition}[theorem]{\indent\sc Proposition}

\theoremstyle{definition} 

\newtheorem{remark}[theorem]{\indent\sc Remark}

\begin{document}

\title[Coupled Painlev\'e III systems]{Coupled Painlev\'e III systems with affine Weyl group symmetry of types $B_4^{(1)}$, $D_4^{(1)}$ and $D_5^{(2)}$ \\}

\author{Yusuke Sasano }

\renewcommand{\thefootnote}{\fnsymbol{footnote}}
\footnote[0]{2000\textit{ Mathematics Subjet Classification}.
14E05,20F55,34M55.}

\maketitle

\begin{abstract}
We find and study four kinds of a 4-parameter family of four-dimensional coupled Painlev\'e III systems with affine Weyl group symmetry of types $B_4^{(1)}$, $D_4^{(1)}$ and $D_5^{(2)}$. We also show that these systems are equivalent by an explicit birational and symplectic transformation, respectively.
\end{abstract}

\section{Introduction}
In \cite{4,5}, we presented some types of coupled Painlev\'e systems with various affine Weyl group symmetries. In this paper, we present a 4-parameter family of 2-coupled Painlev\'e III systems with affine Weyl group symmetry of type $D_4^{(1)}$ explicitly given by
\begin{equation}
\frac{dx}{dt}=\frac{\partial H_{D_4^{(1)}}}{\partial y}, \ \ \frac{dy}{dt}=-\frac{\partial H_{D_4^{(1)}}}{\partial x}, \ \ \frac{dz}{dt}=\frac{\partial H_{D_4^{(1)}}}{\partial w}, \ \ \frac{dw}{dt}=-\frac{\partial H_{D_4^{(1)}}}{\partial z}
\end{equation}
with the Hamiltonian
\begin{align}\label{2}
\begin{split}
H_{D_4^{(1)}} &=H_{III}(x,y,t;\alpha_1,\frac{2\alpha_2+\alpha_3+\alpha_4}{2},\alpha_0)\\
&+{\tilde H}_{III}(z,w,t;\alpha_3,\frac{\alpha_4-\alpha_3}{2},1-\alpha_4)-\frac{2yw}{t}.
\end{split}
\end{align}
Here $x,y,z$ and $w$ denote unknown complex variables and $\alpha_0,\alpha_1,\alpha_2,\alpha_3$ and $\alpha_4$ are complex parameters satisfying the relation $\alpha_0+\alpha_1+2\alpha_2+\alpha_3+\alpha_4=1$. The symbols $H_{III},{\tilde H}_{III}$ are given by
\begin{align}\label{4,5}
&H_{III}(q,p,t;\gamma_0,\gamma_1,\gamma_2)=\frac{q^2p(p-1)+q\{(\gamma_0+\gamma_2)p-\gamma_0\}+tp}{t} \quad (\gamma_0+2\gamma_1+\gamma_2=1),\\
&{\tilde H}_{III}(q,p,t;\gamma_0,\gamma_1,\gamma_2)=\frac{q^2p(p-t)-q\{(-\gamma_0+\gamma_2)p+\gamma_0t\}+p}{t}
\end{align}
with the relation
\begin{align}
\begin{split}\label{6}
&dp \wedge dq - dH_{III}(q,p,t;\gamma_0,\gamma_1,\gamma_2) \wedge dt\\
&=dP \wedge dQ - d{\tilde H}_{III}(Q,P,t;\gamma_0,\gamma_1,\gamma_2) \wedge dt.
\end{split}
\end{align}
Here the relation between $(q,p)$ and $(Q,P)$ is given by
\begin{equation}\label{7}
(Q,P)=(1/q,-q(qp+\gamma_0)).
\end{equation}
We remark that for this system we tried to seek its first integrals of polynomial type with respect to $x,y,z,w$. However, we can not find. Of course, the Hamiltonian $H_{D_4^{(1)}}$ is not the first integral.

The B{\"a}cklund transformations of this system satisfy Noumi-Yamada's universal description for $D_4^{(1)}$ root system (see \cite{2}). Since these universal B{\"a}cklund transformations have Lie theoretic origin, similarity reduction of a Drinfeld-Sokolov hierarchy admits such a B{\"a}cklund symmetry. The aim of this paper is to introduce the system of type $D_4^{(1)}$ and show the relationship between this system and the system of type $B_4^{(1)}$ (see \cite{5}) by an explicit birational and symplectic transformation. We remark that the B{\"a}cklund transformations of that system of type $B_4^{(1)}$ do not have Noumi-Yamada's universal description for $B_4^{(1)}$ root system. In this vein, it had been an open question whether our system of type $B_4^{(1)}$ can be obtained by similarity reduction of a Drinfeld-Sokolov hierarchy. After our discovery of this system, they were studied from the viewpoint of Drinfeld-Sokolov hierarchy by K. Fuji independently (cf. \cite{1}), and he succeeded to obtain our system by similarity reduction of the Drinfeld-Sokolov hierarchy of type $D_4^{(1)}$. His paper will appear soon.
\begin{figure}[h]
\unitlength 0.1in
\begin{picture}(36.83,8.20)(0.27,-19.70)
\put(1.2000,-13.2000){\makebox(0,0)[lb]{Our discovery of the system}}%
\put(1.2000,-15.1000){\makebox(0,0)[lb]{of type $D_4^{(1)}$}}%
%
\put(0.2700,-14.2500){\makebox(0,0)[lb]{}}%
\put(24.1000,-14.4400){\makebox(0,0)[lb]{$W(D_4^{(1)})$}}%
%
\special{pn 20}%
\special{sh 0.600}%
\special{ar 403 1675 12 20  0.0000000 6.2831853}%
\put(4.7300,-17.6200){\makebox(0,0)[lb]{Noumi-Yamada's universal description}}%
\put(37.0000,-17.2000){\makebox(0,0)[lb]{Drinfeld-Sokolov}}%
\put(36.0000,-21.1000){\makebox(0,0)[lb]{(by K. Fuji's work)}}%
\put(4.6700,-19.6100){\makebox(0,0)[lb]{for $D_4^{(1)}$ root system}}%
%
\special{pn 8}%
\special{pa 280 1583}%
\special{pa 280 1970}%
\special{fp}%
\special{pa 280 1550}%
\special{pa 3360 1550}%
\special{fp}%
%
\special{pn 8}%
\special{pa 280 1970}%
\special{pa 3351 1970}%
\special{fp}%
%
\special{pn 8}%
\special{pa 3351 1550}%
\special{pa 3351 1959}%
\special{fp}%
%
\special{pn 20}%
\special{pa 2120 1270}%
\special{pa 2148 1246}%
\special{pa 2176 1225}%
\special{pa 2204 1209}%
\special{pa 2232 1201}%
\special{pa 2261 1202}%
\special{pa 2290 1213}%
\special{pa 2317 1234}%
\special{pa 2340 1260}%
\special{pa 2359 1291}%
\special{pa 2371 1324}%
\special{pa 2377 1358}%
\special{pa 2376 1391}%
\special{pa 2366 1421}%
\special{pa 2346 1446}%
\special{pa 2320 1465}%
\special{pa 2289 1478}%
\special{pa 2257 1483}%
\special{pa 2226 1483}%
\special{pa 2194 1479}%
\special{pa 2162 1473}%
\special{pa 2150 1470}%
\special{sp}%
%
\special{pn 20}%
\special{pa 2190 1480}%
\special{pa 2070 1420}%
\special{fp}%
\special{sh 1}%
\special{pa 2070 1420}%
\special{pa 2121 1468}%
\special{pa 2118 1444}%
\special{pa 2139 1432}%
\special{pa 2070 1420}%
\special{fp}%
%
\special{pn 20}%
\special{pa 3440 1751}%
\special{pa 3680 1751}%
\special{fp}%
\special{sh 1}%
\special{pa 3680 1751}%
\special{pa 3613 1731}%
\special{pa 3627 1751}%
\special{pa 3613 1771}%
\special{pa 3680 1751}%
\special{fp}%
\put(37.1000,-19.0000){\makebox(0,0)[lb]{hierarchy}}%
\end{picture}%
\label{fig:D43}
\caption{}
\end{figure}

Moreover, we presented three kinds of a 4-parameter family of 2-coupled Painlev\'e III systems with extended affine Weyl group symmetry of types $B_4^{(1)}$ and $D_5^{(2)}$ (see \cite{5}), whose Hamiltonians $H_{B_4^{(1)}},\tilde{H}_{B_4^{(1)}}$ and $H_{D_5^{(2)}}$ are given by
\begin{align}
\begin{split}\label{10}
H_{B_4^{(1)}} &={\tilde H}_{III}(x,y,t;\alpha_1,\alpha_2+\frac{\alpha_3+\alpha_4}{2},2\alpha_0+\alpha_1)\\
&+{\tilde H}_{III}(z,w,t;\alpha_3,\frac{\alpha_4-\alpha_3}{2},1-\alpha_4)+\frac{2xw(xy+\alpha_1)}{t},
\end{split}\\
\begin{split}\label{15}
\tilde{H}_{B_4^{(1)}} &=H_{III}(x,y,t;\alpha_1,\alpha_2+\alpha_3+\alpha_4,\alpha_0)\\
&+H_{III}(z,w,t;\alpha_3,\alpha_4,1-\alpha_3-2\alpha_4)+\frac{2yz(zw+\alpha_3)}{t},
\end{split}\\
\begin{split}\label{20}
H_{D_5^{(2)}} &={\tilde H}_{III}(x,y,t;\alpha_1,\alpha_2+\alpha_3+\alpha_4,2\alpha_0+\alpha_1)\\
&+H_{III}(z,w,t;\alpha_3,\alpha_4,1-\alpha_3-2\alpha_4)-\frac{2xz(xy+\alpha_1)(zw+\alpha_3)}{t}.
\end{split}
\end{align}
These systems coincide with the system of type $D_4^{(1)}$ by an explicit birational and symplectic transformation, respectively. In each chart of the phase space, there appear different coupled systems with symmetries of various types.

This paper is organized as follows. In Section 2, we introduce the system of type $D_4^{(1)}$ and its B{\"a}cklund transformations. In Section 3, we introduce two kinds of a 4-parameter family of 2-coupled Painlev\'e III systems with extended affine Weyl group symmetry of type $B_4^{(1)}$ and its B{\"a}cklund transformations. Moreover, these systems coincide with the system of type $D_4^{(1)}$ by an explicit birational and symplectic transformation, respectively. In Section 4, we introduce a 4-parameter family of 2-coupled Painlev\'e III systems with extended affine Weyl group symmetry of type $D_5^{(2)}$ and its B{\"a}cklund transformations. Moreover, this system coincides with the system of type $D_4^{(1)}$ by an explicit birational and symplectic transformation.

\section{The system of type $D_4^{(1)}$}
In this section, we present a 4-parameter family of polynomial Hamiltonian systems that can be considered as 2-coupled Painlev\'e III systems in dimension four given by
\begin{equation}\label{1}
  \left\{
  \begin{aligned}
   \frac{dx}{dt} &=\frac{\partial H_{D_4^{(1)}}}{\partial y}=\frac{2x^2y-{x^2}+(\alpha_0+\alpha_1)x-2w}{t}+1,\\
   \frac{dy}{dt} &=-\frac{\partial H_{D_4^{(1)}}}{\partial x}=\frac{-2xy^2+2xy-(\alpha_0+\alpha_1)y+{\alpha_1}}{t},\\
   \frac{dz}{dt} &=\frac{\partial H_{D_4^{(1)}}}{\partial w}=\frac{2z^2w-{tz^2}-(1-\alpha_3-\alpha_4)z+1-2y}{t},\\
   \frac{dw}{dt} &=-\frac{\partial H_{D_4^{(1)}}}{\partial z}=\frac{-2zw^2+2{tzw}+(1-\alpha_3-\alpha_4)w+\alpha_3t}{t}
   \end{aligned}
  \right. 
\end{equation}
with the Hamiltonian \eqref{2}.
\begin{theorem}\label{2.1}
The system \eqref{1} admits affine Weyl group symmetry of type $D_4^{(1)}$ as the group of its B{\"a}cklund transformations {\rm (cf. \cite{3}), \rm} whose generators are explicitly given as follows{\rm : \rm}with the notation $(*):=(x,y,z,w,t;\alpha_0,\alpha_1,\ldots,\alpha_4)$,
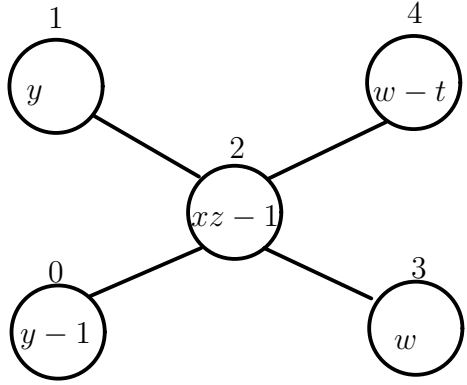
\begin{figure}
\unitlength 0.1in
\begin{picture}(23.71,20.35)(22.18,-23.65)
%
\special{pn 20}%
\special{ar 2462 834 244 244  0.0000000 6.2831853}%
%
\special{pn 20}%
\special{ar 2473 2121 244 244  0.0000000 6.2831853}%
%
\special{pn 20}%
\special{ar 3397 1494 244 244  0.0000000 6.2831853}%
%
\special{pn 20}%
\special{ar 4334 814 244 244  0.0000000 6.2831853}%
%
\special{pn 20}%
\special{ar 4345 2101 244 244  0.0000000 6.2831853}%
%
\special{pn 20}%
\special{pa 2660 988}%
\special{pa 3210 1307}%
\special{fp}%
%
\special{pn 20}%
\special{pa 2638 1934}%
\special{pa 3221 1681}%
\special{fp}%
%
\special{pn 20}%
\special{pa 3573 1318}%
\special{pa 4189 1021}%
\special{fp}%
%
\special{pn 20}%
\special{pa 3551 1681}%
\special{pa 4112 1945}%
\special{fp}%
\put(22.7500,-22.0900){\makebox(0,0)[lb]{$y-1$}}%
\put(23.0800,-9.3300){\makebox(0,0)[lb]{$y$}}%
\put(31.7000,-15.8000){\makebox(0,0)[lb]{$xz-1$}}%
\put(41.2300,-9.2200){\makebox(0,0)[lb]{$w-t$}}%
\put(42.3300,-22.0900){\makebox(0,0)[lb]{$w$}}%
\put(24.2000,-18.6000){\makebox(0,0)[lb]{$0$}}%
\put(24.2000,-5.3000){\makebox(0,0)[lb]{$1$}}%
\put(33.7000,-12.1000){\makebox(0,0)[lb]{$2$}}%
\put(43.2000,-18.4000){\makebox(0,0)[lb]{$3$}}%
\put(43.0000,-5.0000){\makebox(0,0)[lb]{$4$}}%
\end{picture}%
\label{fig:D41}
\caption{Dynkin diagram of type $D_4^{(1)}$}
\end{figure}
\begin{align*}
        s_{0}: (*) &\rightarrow (x+\frac{\alpha_0}{y-1},y,z,w,t;-\alpha_0,\alpha_1,\alpha_2+\alpha_0,\alpha_3,\alpha_4),\\
        s_{1}: (*) &\rightarrow (x+\frac{\alpha_1}{y},y,z,w,t;\alpha_0,-\alpha_1,\alpha_2+\alpha_1,\alpha_3,\alpha_4), \\
        s_{2}: (*) &\rightarrow  (x,y-\frac{\alpha_2z}{xz-1},z,w-\frac{\alpha_2x}{xz-1},t;\alpha_0+\alpha_2,\alpha_1+\alpha_2,-\alpha_2,\alpha_3+\alpha_2,\alpha_4+\alpha_2), \\
        s_{3}: (*) &\rightarrow (x,y,z+\frac{\alpha_3}{w},w,t;\alpha_0,\alpha_1,\alpha_2+\alpha_3,-\alpha_3,\alpha_4), \\
        s_{4}: (*) &\rightarrow (x,y,z+\frac{\alpha_4}{w-t},w,t;\alpha_0,\alpha_1,\alpha_2+\alpha_4,\alpha_3,-\alpha_4), \\
        \pi_1: (*) &\rightarrow (-x,1-y,-z,-w,-t;\alpha_1,\alpha_0,\alpha_2,\alpha_3,\alpha_4), \\
        \pi_2: (*) &\rightarrow (x,y,z,w-t,-t;\alpha_0,\alpha_1,\alpha_2,\alpha_4,\alpha_3), \\
        \pi_3: (*) &\rightarrow (tz,\frac{w}{t},\frac{x}{t},ty,t;\alpha_4,\alpha_3,\alpha_2,\alpha_1,\alpha_0), \\
        \pi_4: (*) &\rightarrow (-tz,\frac{t-w}{t},-\frac{x}{t},t-ty,t;\alpha_3,\alpha_4,\alpha_2,\alpha_0,\alpha_1).
\end{align*}
\end{theorem}

\begin{remark}
The transformations $\pi_2,\pi_3$ and $\pi_4$ satisfy the following relation{\rm : \rm}
\begin{equation}
\pi_4=\pi_2\pi_3\pi_2.
\end{equation}
\end{remark}

\begin{proposition}
Let us define the following translation operators \rm{(see \cite{Ma})\rm}
\begin{align}
\begin{split}
&T_1:=s_3s_0s_2s_4s_1s_2\pi_4, \quad T_2:=s_4s_1s_2s_3s_0s_2\pi_4,\\
&T_3:=s_3s_2s_0s_1s_2s_3\pi_1\pi_2, \quad T_4:=s_4s_3s_2s_1s_0s_2\pi_1\pi_2.
\end{split}
\end{align}
These translation operators act on parameters $\alpha_i$ as follows$:$
\begin{align}
\begin{split}
T_1(\alpha_0,\alpha_1,\ldots,\alpha_4)=&(\alpha_0,\alpha_1,\alpha_2,\alpha_3,\alpha_4)+(1,0,-1,1,0),\\
T_2(\alpha_0,\alpha_1,\ldots,\alpha_4)=&(\alpha_0,\alpha_1,\alpha_2,\alpha_3,\alpha_4)+(0,1,-1,0,1),\\
T_3(\alpha_0,\alpha_1,\ldots,\alpha_4)=&(\alpha_0,\alpha_1,\alpha_2,\alpha_3,\alpha_4)+(0,0,0,1,-1),\\
T_4(\alpha_0,\alpha_1,\ldots,\alpha_4)=&(\alpha_0,\alpha_1,\alpha_2,\alpha_3,\alpha_4)+(0,0,-1,1,1).
\end{split}
\end{align}
\end{proposition}

\begin{theorem}\label{2.2}
Let us consider a polynomial Hamiltonian system with Hamiltonian $H \in {\Bbb C}(t)[x,y,z,w]$. We assume that

$(A1)$ $deg(H)=5$ with respect to $x,y,z,w$.

$(A2)$ This system becomes again a polynomial Hamiltonian system in each coordinate $r_i \ (i=0,1,3,4)${\rm : \rm}
\begin{align*}
&r_0:x_0=1/x, \ y_0=-((y-1)x+\alpha_0)x, \ z_0=z, \ w_0=w, \\
&r_1:x_1=1/x, \ y_1=-(yx+\alpha_1)x, \ z_1=z, \ w_1=w, \\
&r_3:x_3=x, \ y_3=y, \ z_3=1/z, \ w_3=-z(wz+\alpha_3), \\
&r_4:x_4=x, \ y_4=y, \ z_4=1/z, \ w_3=-z((w-t)z+\alpha_4).
\end{align*}

$(A3)$ In addition to the assumption $(A2)$, the Hamiltonian system in the coordinate $r_1$ becomes again a polynomial Hamiltonian system in the coordinate $r_2${\rm : \rm}
\begin{equation*}
r_2:x_2=-((x_1-z_1)y_1-\alpha_2)y_1, \ y_2=1/y_1, \ z_2=z_1, \ w_2=w_1+y_1.
\end{equation*}
Then such a system coincides with the system \eqref{1}.
\end{theorem}
Each coordinate $r_i \ (i=0,1,3,4)$ contains a three-parameter family of meromorphic solutions of \eqref{1}.

Theorems \ref{2.1} and \ref{2.2} can be checked by a direct calculation, respectively.

We note that the following transformations
\begin{align*}
        w_{0}: (*) &\rightarrow (x+\frac{\alpha_0}{y-1},y,z,w,t;-\alpha_0,\alpha_1,\alpha_2+\alpha_0,\alpha_3,\alpha_4), \\
        w_{1}: (*) &\rightarrow (x+\frac{\alpha_1}{y},y,z,w,t;\alpha_0,-\alpha_1,\alpha_2+\alpha_1,\alpha_3,\alpha_4), \\
        w_{2}: (*) &\rightarrow  (x,y-\frac{\alpha_2}{x-z},z,w+\frac{\alpha_2}{x-z},t;\alpha_0+\alpha_2,\alpha_1+\alpha_2,-\alpha_2,\alpha_3+\alpha_2,\alpha_4+\alpha_2), \\
        w_{3}: (*) &\rightarrow (x,y,z+\frac{\alpha_3}{w},w,t;\alpha_0,\alpha_1,\alpha_2+\alpha_3,-\alpha_3,\alpha_4), \\
        w_{4}: (*) &\rightarrow (x,y,z+\frac{\alpha_4}{w-t},w,t;\alpha_0,\alpha_1,\alpha_2+\alpha_4,\alpha_3,-\alpha_4)
\end{align*}
define a representation of the affine Weyl group of type $D_4^{(1)}$. However, we can not find polynomial Hamiltonian systems with affine Weyl group symmetry of type $D_4^{(1)}$ described above.

Moreover, from the viewpoint of holomorphy conditions let us consider a polynomial Hamiltonian system with $H \in {\Bbb C}(t)[x,y,z,w]$. We assume that

$(A)$ This system becomes again a polynomial Hamiltonian system in each coordinate $r_i \ (i=0,1, \dots ,4)${\rm : \rm}
\begin{align*}
&r_0:x_0=1/x, \ y_0=-((y-1)x+\alpha_0)x, \ z_0=z, \ w_0=w, \\
&r_1:x_1=1/x, \ y_1=-(yx+\alpha_1)x, \ z_1=z, \ w_1=w, \\
&r_2:x_2=-((x-z)y-\alpha_2)y, \ y_2=1/y, \ z_2=z, \ w_2=w+y,\\
&r_3:x_3=x, \ y_3=y, \ z_3=1/z, \ w_3=-z(wz+\alpha_3), \\
&r_4:x_4=x, \ y_4=y, \ z_4=1/z, \ w_3=-z((w-t)z+\alpha_4).
\end{align*}
It is still an open question whether we can find a system satisfying the assumption $(A)$.

We also give an explicit description of a confluence from 2-coupled Painlev\'e V system with $W(D_5^{(1)})$-symmetry to the system of type $D_4^{(1)}$. At first, we recall 5-parameter family of 2-coupled Painlev\'e V systems with $W(D_5^{(1)})$-symmetry (see \cite{5}) explicitly given by
\begin{equation}\label{D5}
\frac{dx}{dt}=\frac{\partial H_{D_5^{(1)}}}{\partial y}, \ \ \frac{dy}{dt}=-\frac{\partial H_{D_5^{(1)}}}{\partial x}, \ \ \frac{dz}{dt}=\frac{\partial H_{D_5^{(1)}}}{\partial w}, \ \ \frac{dw}{dt}=-\frac{\partial H_{D_5^{(1)}}}{\partial z}
\end{equation}
with the Hamiltonian
\begin{align}\label{HD5}
\begin{split}
H_{D_5^{(1)}}=&H_{V}(x,y,t;\beta_2+\beta_5,\beta_1,\beta_2+2\beta_3+\beta_4)+H_{V}(z,w,t;\beta_5,\beta_3,\beta_4)\\
&+\frac{2yz\{(z-1)w+\beta_3\}}{t},
\end{split}
\end{align}
where the symbol $H_{V}(q,p,t;\gamma_1,\gamma_2,\gamma_3)$ denotes the Hamiltonian of the second-order Painlev\'e V systems given by
$$
H_{V}(q,p,t;\gamma_1,\gamma_2,\gamma_3)=\frac{q(q-1)p(p+t)-(\gamma_1+\gamma_3)qp+\gamma_1p+\gamma_2tq}{t}.
$$
Here $\beta_0,\beta_1,\dots,\beta_5$ are complex parameters normalized as $\beta_0+\beta_1+2\beta_2+2\beta_3+\beta_4+\beta_5=1$.

The system \eqref{D5} admits affine Weyl group symmetry of type $D_5^{(1)}$ as the group of its B{\"a}cklund transformations, whose generators $w_0,w_1,\dots,w_5$ defined as follows$:$ with {\it the notation} $(*):=(x,y,z,w,t;\beta_0,\beta_1,\dots,\beta_5)$,
\begin{align}
\begin{split}
w_0: (*) \rightarrow &(x+\frac{\beta_0}{y+t},y,z,w,t;-\beta_0,\beta_1,\beta_2+\beta_0,\beta_3,\beta_4,\beta_5),\\
w_1: (*) \rightarrow &(x+\frac{\beta_1}{y},y,z,w,t;\beta_0,-\beta_1,\beta_2+\beta_1,\beta_3,\beta_4,\beta_5),\\
w_2: (*) \rightarrow &(x,y-\frac{\beta_2}{x-z},z,w+\frac{\beta_2}{x-z},t;\beta_0+\beta_2,\beta_1+\beta_2,-\beta_2,\beta_3+\beta_2,\beta_4,\beta_5),\\
w_3: (*) \rightarrow &(x,y,z+\frac{\beta_3}{w},w,t;\beta_0,\beta_1,\beta_2+\beta_3,-\beta_3,\beta_4+\beta_3,\beta_5+\beta_3),\\
w_4: (*) \rightarrow &(x,y,z,w-\frac{\beta_4}{(z-1)},t;\beta_0,\beta_1,\beta_2,\beta_3+\beta_4,-\beta_4,\beta_5),\\
w_5: (*) \rightarrow &(x,y,z,w-\frac{\beta_5}{z},t;\beta_0,\beta_1,\beta_2,\beta_3+\beta_5,\beta_4,-\beta_5).
\end{split}
\end{align}

\begin{proposition}\label{th:1.3}
For the system of type $D_5^{(1)}$, we make the change of parameters and variables
\begin{gather}
\begin{gathered}
\beta_0=\alpha_0, \ \beta_1=\alpha_1, \ \beta_2=\alpha_2, \ \beta_3=\alpha_3, \ \beta_4=\alpha_4-\alpha_3-\frac{1}{\varepsilon}, \ \beta_5=\frac{1}{\varepsilon},
\end{gathered}\\
\begin{gathered}
t=-\varepsilon T, \ x=1+\frac{X}{\varepsilon T}, \ y=\varepsilon TY, \ z=1+\frac{1}{\varepsilon TZ}, \ w=-\varepsilon T(ZW+A_3)Z
\end{gathered}
\end{gather}
from $\beta_0,\beta_1,\dots,\beta_5,t,x,y,z,w$ to $\alpha_0,\alpha_1,\dots,\alpha_4,\varepsilon,T,X,Y,Z,W$. Then the system can also be written in the new variables $T,X,Y,Z,W$ and parameters $\alpha_0,\alpha_1,\ldots,\alpha_4,\varepsilon$ as a Hamiltonian system. This new system tends to the system \eqref{1} of type $D_4^{(1)}$ as $\varepsilon \rightarrow 0$.
\end{proposition}

By proving the following theorem, we see how the degeneration process in Proposition \ref{th:1.3} works on the B{\"a}cklund transformation group $W(D_5^{(1)})=<w_0,w_1,\ldots,w_5>$ described above.
\begin{proposition}
For the degeneration process in Proposition \ref{th:1.3}, we can choose a subgroup 
$$
W_{D_5^{(1)} \rightarrow D_4^{(1)}}:=\{<s_0,\ldots,s_4>|s_i:=w_i \ (i=0,1,2,3), \ s_4:=w_4w_5w_3w_4w_5\}
$$
of the B{\"a}cklund transformation group $W(D_5^{(1)})$ so that $W_{D_5^{(1)} \rightarrow D_4^{(1)}}$ converges to $W(D_4^{(1)})$ as $\varepsilon \rightarrow 0$.
\end{proposition}

\section{The system of type $B_4^{(1)}$ }
In this section, we propose two types of a 4-parameter family of 2-coupled Painlev\'e III systems in dimension four with affine Weyl group symmetry of type $B_4^{(1)}$. Each of them is equivalent to a polynomial Hamiltonian system, however, each has a different representaion of type $B_4^{(1)}$. We also show that each of them is equivalent to the system \eqref{1} by a birational and symplectic transformation.

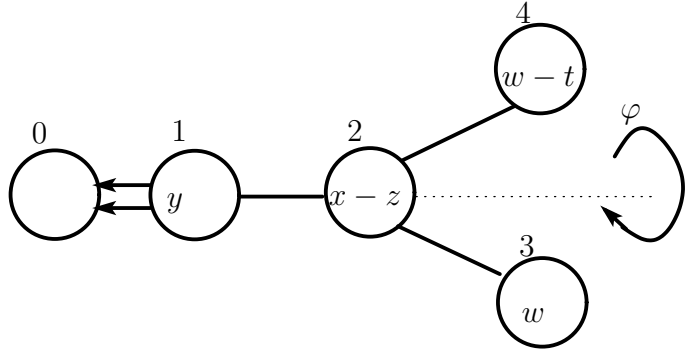
\begin{figure}[h]
\unitlength 0.1in
\begin{picture}(35.20,18.66)(16.38,-23.56)
%
\special{pn 20}%
\special{ar 2602 1562 232 232  0.0000000 6.2831853}%
%
\special{pn 20}%
\special{ar 1870 1560 232 232  0.0000000 6.2831853}%
%
\special{pn 20}%
\special{ar 3520 1548 232 232  0.0000000 6.2831853}%
%
\special{pn 20}%
\special{ar 4410 902 232 232  0.0000000 6.2831853}%
%
\special{pn 20}%
\special{ar 4421 2124 232 232  0.0000000 6.2831853}%
%
\special{pn 20}%
\special{pa 3687 1381}%
\special{pa 4272 1098}%
\special{fp}%
%
\special{pn 20}%
\special{pa 3666 1725}%
\special{pa 4199 1976}%
\special{fp}%
\put(24.5600,-16.5600){\makebox(0,0)[lb]{$y$}}%
\put(33.0400,-16.3000){\makebox(0,0)[lb]{$x-z$}}%
\put(42.1000,-10.0400){\makebox(0,0)[lb]{$w-t$}}%
\put(43.1400,-22.2700){\makebox(0,0)[lb]{$w$}}%
\put(17.5000,-12.9000){\makebox(0,0)[lb]{$0$}}%
\put(24.8000,-12.8000){\makebox(0,0)[lb]{$1$}}%
\put(34.0000,-12.8000){\makebox(0,0)[lb]{$2$}}%
\put(43.0000,-18.8000){\makebox(0,0)[lb]{$3$}}%
\put(42.8000,-6.6000){\makebox(0,0)[lb]{$4$}}%
%
\special{pn 20}%
\special{pa 2840 1570}%
\special{pa 3270 1570}%
\special{fp}%
%
\special{pn 20}%
\special{pa 2370 1510}%
\special{pa 2100 1510}%
\special{fp}%
\special{sh 1}%
\special{pa 2100 1510}%
\special{pa 2167 1530}%
\special{pa 2153 1510}%
\special{pa 2167 1490}%
\special{pa 2100 1510}%
\special{fp}%
%
\special{pn 20}%
\special{pa 2370 1630}%
\special{pa 2100 1630}%
\special{fp}%
\special{sh 1}%
\special{pa 2100 1630}%
\special{pa 2167 1650}%
\special{pa 2153 1630}%
\special{pa 2167 1610}%
\special{pa 2100 1630}%
\special{fp}%
%
\special{pn 8}%
\special{pa 3760 1570}%
\special{pa 4990 1570}%
\special{dt 0.045}%
\special{pa 4990 1570}%
\special{pa 4989 1570}%
\special{dt 0.045}%
%
\special{pn 20}%
\special{pa 4800 1360}%
\special{pa 4820 1326}%
\special{pa 4841 1294}%
\special{pa 4862 1264}%
\special{pa 4883 1240}%
\special{pa 4906 1222}%
\special{pa 4930 1212}%
\special{pa 4955 1211}%
\special{pa 4981 1219}%
\special{pa 5008 1235}%
\special{pa 5034 1258}%
\special{pa 5059 1284}%
\special{pa 5082 1314}%
\special{pa 5103 1345}%
\special{pa 5120 1375}%
\special{pa 5134 1406}%
\special{pa 5145 1437}%
\special{pa 5153 1468}%
\special{pa 5157 1500}%
\special{pa 5158 1532}%
\special{pa 5155 1564}%
\special{pa 5149 1597}%
\special{pa 5140 1630}%
\special{pa 5127 1664}%
\special{pa 5111 1697}%
\special{pa 5093 1728}%
\special{pa 5072 1755}%
\special{pa 5049 1777}%
\special{pa 5023 1793}%
\special{pa 4997 1800}%
\special{pa 4969 1799}%
\special{pa 4940 1791}%
\special{pa 4911 1776}%
\special{pa 4881 1758}%
\special{pa 4850 1737}%
\special{pa 4840 1730}%
\special{sp}%
%
\special{pn 20}%
\special{pa 4860 1740}%
\special{pa 4750 1640}%
\special{fp}%
\special{sh 1}%
\special{pa 4750 1640}%
\special{pa 4786 1700}%
\special{pa 4789 1676}%
\special{pa 4813 1670}%
\special{pa 4750 1640}%
\special{fp}%
\put(48.3000,-11.9000){\makebox(0,0)[lb]{$\varphi$}}%
\end{picture}%
\label{fig:D44}
\caption{Dynkin diagram of type $B_4^{(1)}$}
\end{figure}
The first member is given by
\begin{equation}\label{9}
  \left\{
  \begin{aligned}
   \frac{dx}{dt} &=\frac{\partial H_{B_4^{(1)}}}{\partial y}=\frac{2x^2y-tx^2-2\alpha_0x+1}{t}+\frac{2x^2w}{t},\\
   \frac{dy}{dt} &=-\frac{\partial H_{B_4^{(1)}}}{\partial x}=\frac{-2xy^2+2txy+2\alpha_0y+\alpha_1t}{t}-\frac{2w(2xy+\alpha_1)}{t},\\
   \frac{dz}{dt} &=\frac{\partial H_{B_4^{(1)}}}{\partial w}=\frac{2z^2w-tz^2-(1-\alpha_3-\alpha_4)z+1}{t}+\frac{2x(xy+\alpha_1)}{t},\\
   \frac{dw}{dt} &=-\frac{\partial H_{B_4^{(1)}}}{\partial z}=\frac{-2zw^2+2tzw+(1-\alpha_3-\alpha_4)w+\alpha_3t}{t}
   \end{aligned}
  \right. 
\end{equation}
with the Hamiltonian \eqref{10}. Here $x,y,z$ and $w$ denote unknown complex variables and $\alpha_0,\alpha_1,\alpha_2,\alpha_3$ and $\alpha_4$ are complex parameters satisfying the relation $2\alpha_0+2\alpha_1+2\alpha_2+\alpha_3+\alpha_4=1$.

\begin{theorem}\label{3.1}
The system \eqref{9} admits extended affine Weyl group symmetry of type $B_4^{(1)}$ as the group of its B{\"a}cklund transformations {\rm (cf. \cite{3}), \rm} whose generators are explicitly given as follows{\rm : \rm}with the notation $(*):=(x,y,z,w,t;\alpha_0,\alpha_1,\ldots,\alpha_4)$,
\begin{align*}
        s_{0}: (*) &\rightarrow (-x,-y+\frac{2\alpha_0}{x}-\frac{1}{x^2},-z,-w,-t;-\alpha_0,\alpha_1+2\alpha_0,\alpha_2,\alpha_3,\alpha_4),\\
        s_{1}: (*) &\rightarrow (x+\frac{\alpha_1}{y},y,z,w,t;\alpha_0+\alpha_1,-\alpha_1,\alpha_2+\alpha_1,\alpha_3,\alpha_4), \\
        s_{2}: (*) &\rightarrow  (x,y-\frac{\alpha_2}{x-z},z,w+\frac{\alpha_2}{x-z},t;\alpha_0,\alpha_1+\alpha_2,-\alpha_2,\alpha_3+\alpha_2,\alpha_4+\alpha_2), \\
        s_{3}: (*) &\rightarrow (x,y,z+\frac{\alpha_3}{w},w,t;\alpha_0,\alpha_1,\alpha_2+\alpha_3,-\alpha_3,\alpha_4), \\
        s_{4}: (*) &\rightarrow (x,y,z+\frac{\alpha_4}{w-t},w,t;\alpha_0,\alpha_1,\alpha_2+\alpha_4,\alpha_3,-\alpha_4), \\
        \varphi: (*) &\rightarrow (x,y,z,w-t,-t;\alpha_0,\alpha_1,\alpha_2,\alpha_4,\alpha_3).
\end{align*}
\end{theorem}

\begin{theorem}\label{3.2}
Let us consider a polynomial Hamiltonian system with Hamiltonian $H \in {\Bbb C}(t)[x,y,z,w]$. We assume that

$(A1)$ $deg(H)=5$ with respect to $x,y,z,w$.

$(A2)$ This system becomes again a polynomial Hamiltonian system in each coordinate $r_i \ (i=0,1,\dots,4)${\rm : \rm}
\begin{align*}
&r_0:x_0=x, \ y_0=y-\frac{2\alpha_0}{x}+\frac{1}{x^2}, \ z_0=z, \ w_0=w, \\
&r_1:x_1=1/x, \ y_1=-(yx+\alpha_1)x, \ z_1=z, \ w_1=w, \\
&r_2:x_2=-((x-z)y-\alpha_2)y, \ y_2=1/y, \ z_2=z, \ w_2=w+y,\\
&r_3:x_3=x, \ y_3=y, \ z_3=1/z, \ w_3=-z(wz+\alpha_3), \\
&r_4:x_4=x, \ y_4=y, \ z_4=1/z, \ w_4=-z((w-t)z+\alpha_4).
\end{align*}
Then such a system coincides with the system \eqref{9}.
\end{theorem}

Theorems \ref{3.1} and \ref{3.2} can be checked by a direct calculation, respectively.

\begin{theorem}\label{3.3}
For the system \eqref{1} of type $D_4^{(1)}$, we make the change of parameters and variables
\begin{gather}
\begin{gathered}\label{12}
A_0=\frac{\alpha_0-\alpha_1}{2}, \quad A_1=\alpha_1, \quad A_2=\alpha_2, \quad A_3=\alpha_3, \quad A_4=\alpha_4,\\
\end{gathered}\\
\begin{gathered}\label{13}
X=\frac{1}{x}, \quad Y=-(xy+\alpha_1)x, \quad Z=z, \quad W=w
\end{gathered}
\end{gather}
from $\alpha_0,\alpha_1, \dots ,\alpha_4,x,y,z,w$ to $A_0,A_1,\dots ,A_4,X,Y,Z,W$. Then the system \eqref{1} can also be written in the new variables $X,Y,Z,W$ and parameters $A_0,A_1,\dots ,A_4$ as a Hamiltonian system. This new system tends to the system \eqref{9} with the Hamiltonian \eqref{10}.
\end{theorem}

\begin{proof}
Notice that
$$
2A_0+2A_1+2A_2+A_3+A_4=\alpha_0+\alpha_1+2\alpha_2+\alpha_3+\alpha_4=1
$$
and the change of variables from $(x,y,z,w)$ to $(X,Y,Z,W)$ is symplectic. Choose $S_i \ (i=0,1,\dots ,4)$ and $\varphi$ as
$$
S_0:=\pi_1, \ S_1:=s_1, \ S_2:=s_2, \ S_3:=s_3, \ S_4:=s_4, \ \varphi:=\pi_2.
$$
Then the transformations $S_i$ are reflections of the parameters $A_0,A_1,\dots ,A_4$. The transformation group $\tilde{W}(B_4^{(1)})=<S_0,S_1,\dots ,S_4,\varphi>$ coincides with the transformations given in Theorem \ref{3.1}.
\end{proof}

The second member (see (4),(5)) is given by
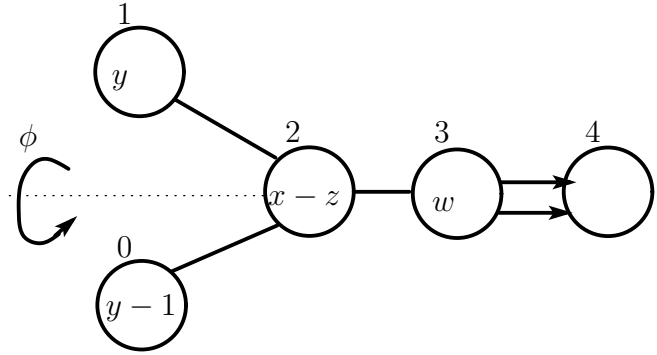
\begin{figure}[h]
\unitlength 0.1in
\begin{picture}(33.66,18.76)(15.30,-24.16)
%
\special{pn 20}%
\special{ar 2212 962 232 232  0.0000000 6.2831853}%
%
\special{pn 20}%
\special{ar 2223 2184 232 232  0.0000000 6.2831853}%
%
\special{pn 20}%
\special{ar 3101 1589 232 232  0.0000000 6.2831853}%
%
\special{pn 20}%
\special{pa 2400 1108}%
\special{pa 2923 1411}%
\special{fp}%
%
\special{pn 20}%
\special{pa 2380 2007}%
\special{pa 2933 1766}%
\special{fp}%
\put(20.3500,-22.6800){\makebox(0,0)[lb]{$y-1$}}%
\put(20.6600,-10.5600){\makebox(0,0)[lb]{$y$}}%
\put(28.8500,-16.7000){\makebox(0,0)[lb]{$x-z$}}%
%
\special{pn 20}%
\special{pa 3344 1590}%
\special{pa 3624 1590}%
\special{fp}%
%
\special{pn 20}%
\special{ar 3874 1600 232 232  0.0000000 6.2831853}%
%
\special{pn 20}%
\special{ar 4664 1590 232 232  0.0000000 6.2831853}%
%
\special{pn 20}%
\special{pa 4104 1540}%
\special{pa 4464 1540}%
\special{fp}%
\special{sh 1}%
\special{pa 4464 1540}%
\special{pa 4397 1520}%
\special{pa 4411 1540}%
\special{pa 4397 1560}%
\special{pa 4464 1540}%
\special{fp}%
%
\special{pn 20}%
\special{pa 4094 1700}%
\special{pa 4444 1700}%
\special{fp}%
\special{sh 1}%
\special{pa 4444 1700}%
\special{pa 4377 1680}%
\special{pa 4391 1700}%
\special{pa 4377 1720}%
\special{pa 4444 1700}%
\special{fp}%
\put(37.4400,-16.9000){\makebox(0,0)[lb]{$w$}}%
\put(20.9400,-19.3000){\makebox(0,0)[lb]{$0$}}%
\put(20.9400,-7.1000){\makebox(0,0)[lb]{$1$}}%
\put(29.7400,-13.3000){\makebox(0,0)[lb]{$2$}}%
\put(37.5400,-13.3000){\makebox(0,0)[lb]{$3$}}%
\put(45.4400,-13.3000){\makebox(0,0)[lb]{$4$}}%
%
\special{pn 8}%
\special{pa 1530 1610}%
\special{pa 2860 1610}%
\special{dt 0.045}%
\special{pa 2860 1610}%
\special{pa 2859 1610}%
\special{dt 0.045}%
%
\special{pn 20}%
\special{pa 1840 1470}%
\special{pa 1809 1451}%
\special{pa 1777 1434}%
\special{pa 1747 1421}%
\special{pa 1717 1416}%
\special{pa 1688 1421}%
\special{pa 1661 1435}%
\special{pa 1637 1458}%
\special{pa 1617 1485}%
\special{pa 1602 1513}%
\special{pa 1591 1543}%
\special{pa 1584 1574}%
\special{pa 1581 1605}%
\special{pa 1579 1638}%
\special{pa 1579 1672}%
\special{pa 1580 1706}%
\special{pa 1581 1741}%
\special{pa 1585 1774}%
\special{pa 1594 1805}%
\special{pa 1610 1830}%
\special{pa 1637 1849}%
\special{pa 1669 1859}%
\special{pa 1702 1860}%
\special{pa 1732 1850}%
\special{pa 1759 1834}%
\special{pa 1786 1814}%
\special{pa 1790 1810}%
\special{sp}%
%
\special{pn 20}%
\special{pa 1770 1820}%
\special{pa 1860 1740}%
\special{fp}%
\special{sh 1}%
\special{pa 1860 1740}%
\special{pa 1797 1769}%
\special{pa 1820 1775}%
\special{pa 1823 1799}%
\special{pa 1860 1740}%
\special{fp}%
\put(15.8000,-13.8000){\makebox(0,0)[lb]{$\phi$}}%
\end{picture}%
\label{fig:D45}
\caption{Dynkin diagram of type $B_4^{(1)}$}
\end{figure}
\begin{equation}\label{14}
  \left\{
  \begin{aligned}
   \frac{dx}{dt} &=\frac{\partial \tilde{H}_{B_4^{(1)}}}{\partial y}=\frac{2x^2y-x^2+(\alpha_0+\alpha_1)x+t}{t}+\frac{2z(zw+\alpha_3)}{t},\\
   \frac{dy}{dt} &=-\frac{\partial \tilde{H}_{B_4^{(1)}}}{\partial x}=\frac{-2xy^2+2xy-(\alpha_0+\alpha_1)y+\alpha_1}{t},\\
   \frac{dz}{dt} &=\frac{\partial \tilde{H}_{B_4^{(1)}}}{\partial w}=\frac{2z^2w-z^2+(1-2\alpha_4)z+t}{t}+\frac{2yz^2}{t},\\
   \frac{dw}{dt} &=-\frac{\partial \tilde{H}_{B_4^{(1)}}}{\partial z}=\frac{-2zw^2+2zw-(1-2\alpha_4)w+\alpha_3}{t}-\frac{2y(2zw+\alpha_3)}{t}
   \end{aligned}
  \right. 
\end{equation}
with the Hamiltonian \eqref{15}. Here $x,y,z$ and $w$ denote unknown complex variables and $\alpha_0,\alpha_1,\ldots,\alpha_4$ are complex parameters satisfying the relation $\alpha_0+\alpha_1+2\alpha_2+2\alpha_3+2\alpha_4=1$.

\begin{theorem}\label{3.4}
The system \eqref{14} admits extended affine Weyl group symmetry of type $B_4^{(1)}$ as the group of its B{\"a}cklund transformations {\rm (cf. \cite{3}), \rm} whose generators are explicitly given as follows{\rm : \rm}with the notation $(*):=(x,y,z,w,t;\alpha_0,\alpha_1,\ldots,\alpha_4)$,
\begin{align*}
        s_{0}: (*) &\rightarrow (x+\frac{\alpha_0}{y-1},y,z,w,t;-\alpha_0,\alpha_1,\alpha_2+\alpha_0,\alpha_3,\alpha_4),\\
        s_{1}: (*) &\rightarrow (x+\frac{\alpha_1}{y},y,z,w,t;\alpha_0,-\alpha_1,\alpha_2+\alpha_1,\alpha_3,\alpha_4), \\
        s_{2}: (*) &\rightarrow  (x,y-\frac{\alpha_2}{x-z},z,w+\frac{\alpha_2}{x-z},t;\alpha_0+\alpha_2,\alpha_1+\alpha_2,-\alpha_2,\alpha_3+\alpha_2,\alpha_4), \\
        s_{3}: (*) &\rightarrow (x,y,z+\frac{\alpha_3}{w},w,t;\alpha_0,\alpha_1,\alpha_2+\alpha_3,-\alpha_3,\alpha_4+\alpha_3), \\
        s_{4}: (*) &\rightarrow (x,y,z,w-\frac{2\alpha_4}{z}+\frac{t}{z^2},-t;\alpha_0,\alpha_1,\alpha_2,\alpha_3+2\alpha_4,-\alpha_4), \\
        \phi: (*) &\rightarrow (-x,1-y,-z,-w,-t;\alpha_1,\alpha_0,\alpha_2,\alpha_3,\alpha_4).
\end{align*}
\end{theorem}

\begin{theorem}\label{3.5}
Let us consider a polynomial Hamiltonian system with Hamiltonian $H \in {\Bbb C}(t)[x,y,z,w]$. We assume that

$(A1)$ $deg(H)=5$ with respect to $x,y,z,w$.

$(A2)$ This system becomes again a polynomial Hamiltonian system in each coordinate $r_i \ (i=0,1,\dots,4)${\rm : \rm}
\begin{align*}
&r_0:x_0=1/x, \ y_0=-((y-1)x+\alpha_0)x, \ z_0=z, \ w_0=w, \\
&r_1:x_1=1/x, \ y_1=-(yx+\alpha_1)x, \ z_1=z, \ w_1=w, \\
&r_2:x_2=-((x-z)y-\alpha_2)y, \ y_2=1/y, \ z_2=z, \ w_2=w+y,\\
&r_3:x_3=x, \ y_3=y, \ z_3=1/z, \ w_3=-z(wz+\alpha_3), \\
&r_4:x_4=x, \ y_4=y, \ z_4=z, \ w_4=w-\frac{2\alpha_4}{z}+\frac{t}{z^2}.
\end{align*}
Then such a system coincides with the system \eqref{14}.
\end{theorem}

Theorems \ref{3.4} and \ref{3.5} can be checked by a direct calculation, respectively.

\begin{theorem}\label{3.6}
For the system \eqref{1} of type $D_4^{(1)}$, we make the change of parameters and variables
\begin{gather}
\begin{gathered}\label{17}
A_0=\alpha_0, \quad A_1=\alpha_1, \quad A_2=\alpha_2, \quad A_3=\alpha_3, \quad A_4=\frac{\alpha_4-\alpha_3}{2},\\
\end{gathered}\\
\begin{gathered}\label{18}
X=x, \quad Y=y, \quad Z=\frac{1}{z}, \quad W=-(zw+\alpha_3)z
\end{gathered}
\end{gather}
from $\alpha_0,\alpha_1, \dots ,\alpha_4,x,y,z,w$ to $A_0,A_1,\dots ,A_4,X,Y,Z,W$. Then the system \eqref{1} can also be written in the new variables $X,Y,Z,W$ and parameters $A_0,A_1,\dots ,A_4$ as a Hamiltonian system. This new system tends to the system \eqref{14} with the Hamiltonian \eqref{15}.
\end{theorem}

\begin{proof}
Notice that
$$
A_0+A_1+2A_2+2A_3+2A_4=\alpha_0+\alpha_1+2\alpha_2+\alpha_3+\alpha_4=1
$$
and the change of variables from $(x,y,z,w)$ to $(X,Y,Z,W)$ is symplectic. Choose $S_i \ (i=0,1,\dots ,4)$ and $\phi$ as
$$
S_0:=s_0, \ S_1:=s_1, \ S_2:=s_2, \ S_3:=s_3, \ S_4:=\pi_1, \ \phi:=\pi_2.
$$
Then the transformations $S_i$ are reflections of the parameters $A_0,A_1,\dots ,A_4$. The transformation group $\tilde{W}(B_4^{(1)})=<S_0,S_1,\dots ,S_4,\phi>$ coincides with the transformations given in Theorem \ref{3.4}.
\end{proof}

By using Theorems \ref{3.3} and \ref{3.6}, it is easy to see that the system \eqref{9} coincides with the system \eqref{14} by an explicit birational and symplectic transformation.
\begin{proposition}
For the system \eqref{9} of type $B_4^{(1)}$, we make the change of parameters and variables
\begin{gather*}
\begin{gathered}\label{pro1}
A_0=2\alpha_0+\alpha_1, \quad A_1=\alpha_1, \quad A_2=\alpha_2, \quad A_3=\alpha_3, \quad A_4=\frac{\alpha_4-\alpha_3}{2},\\
\end{gathered}\\
\begin{gathered}\label{pro2}
X=\frac{1}{x}, \quad Y=-(xy+\alpha_1)x, \quad Z=\frac{1}{z}, \quad W=-(zw+\alpha_3)z
\end{gathered}
\end{gather*}
from $\alpha_0,\alpha_1, \dots ,\alpha_4,x,y,z,w$ to $A_0,A_1,\dots ,A_4,X,Y,Z,W$. Then the system \eqref{9} can also be written in the new variables $X,Y,Z,W$ and parameters $A_0,A_1,\dots ,A_4$ as a Hamiltonian system. This new system tends to the system \eqref{14} with the Hamiltonian \eqref{15}.
\end{proposition}

\section{The system of type $D_5^{(2)}$ }
In this section, we propose a 4-parameter family of 2-coupled Painlev\'e III systems in dimension four with affine Weyl group symmetry of type $D_5^{(2)}$ given by
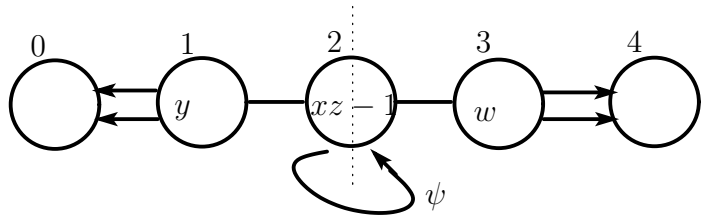
\begin{figure}[ht]
\unitlength 0.1in
\begin{picture}(36.02,10.82)(16.10,-17.62)
%
\special{pn 20}%
\special{ar 2612 1182 232 232  0.0000000 6.2831853}%
%
\special{pn 20}%
\special{ar 1842 1194 232 232  0.0000000 6.2831853}%
%
\special{pn 20}%
\special{ar 3392 1179 232 232  0.0000000 6.2831853}%
\put(24.6600,-12.7600){\makebox(0,0)[lb]{$y$}}%
\put(31.8000,-12.6000){\makebox(0,0)[lb]{$xz-1$}}%
%
\special{pn 20}%
\special{pa 3635 1180}%
\special{pa 3915 1180}%
\special{fp}%
%
\special{pn 20}%
\special{ar 4165 1190 232 232  0.0000000 6.2831853}%
%
\special{pn 20}%
\special{ar 4980 1180 232 232  0.0000000 6.2831853}%
%
\special{pn 20}%
\special{pa 4395 1130}%
\special{pa 4755 1130}%
\special{fp}%
\special{sh 1}%
\special{pa 4755 1130}%
\special{pa 4688 1110}%
\special{pa 4702 1130}%
\special{pa 4688 1150}%
\special{pa 4755 1130}%
\special{fp}%
\put(40.3500,-12.8000){\makebox(0,0)[lb]{$w$}}%
\put(17.1300,-9.4000){\makebox(0,0)[lb]{$0$}}%
\put(24.9400,-9.3000){\makebox(0,0)[lb]{$1$}}%
\put(32.6500,-9.2000){\makebox(0,0)[lb]{$2$}}%
\put(40.4500,-9.2000){\makebox(0,0)[lb]{$3$}}%
\put(48.3500,-9.2000){\makebox(0,0)[lb]{$4$}}%
%
\special{pn 20}%
\special{pa 2860 1180}%
\special{pa 3140 1180}%
\special{fp}%
%
\special{pn 20}%
\special{pa 2370 1120}%
\special{pa 2070 1120}%
\special{fp}%
\special{sh 1}%
\special{pa 2070 1120}%
\special{pa 2137 1140}%
\special{pa 2123 1120}%
\special{pa 2137 1100}%
\special{pa 2070 1120}%
\special{fp}%
%
\special{pn 20}%
\special{pa 2380 1270}%
\special{pa 2080 1270}%
\special{fp}%
\special{sh 1}%
\special{pa 2080 1270}%
\special{pa 2147 1290}%
\special{pa 2133 1270}%
\special{pa 2147 1250}%
\special{pa 2080 1270}%
\special{fp}%
%
\special{pn 20}%
\special{pa 4400 1270}%
\special{pa 4760 1270}%
\special{fp}%
\special{sh 1}%
\special{pa 4760 1270}%
\special{pa 4693 1250}%
\special{pa 4707 1270}%
\special{pa 4693 1290}%
\special{pa 4760 1270}%
\special{fp}%
%
\special{pn 20}%
\special{pa 3266 1440}%
\special{pa 3226 1447}%
\special{pa 3189 1455}%
\special{pa 3155 1464}%
\special{pa 3126 1475}%
\special{pa 3105 1489}%
\special{pa 3093 1507}%
\special{pa 3090 1528}%
\special{pa 3097 1552}%
\special{pa 3112 1578}%
\special{pa 3134 1605}%
\special{pa 3162 1632}%
\special{pa 3195 1658}%
\special{pa 3232 1682}%
\special{pa 3273 1703}%
\special{pa 3316 1721}%
\special{pa 3360 1736}%
\special{pa 3405 1747}%
\special{pa 3450 1755}%
\special{pa 3494 1760}%
\special{pa 3536 1762}%
\special{pa 3575 1761}%
\special{pa 3612 1757}%
\special{pa 3644 1750}%
\special{pa 3671 1741}%
\special{pa 3693 1729}%
\special{pa 3708 1714}%
\special{pa 3716 1697}%
\special{pa 3716 1678}%
\special{pa 3707 1657}%
\special{pa 3691 1634}%
\special{pa 3670 1609}%
\special{pa 3644 1584}%
\special{pa 3616 1558}%
\special{pa 3596 1540}%
\special{sp}%
%
\special{pn 20}%
\special{pa 3626 1560}%
\special{pa 3516 1450}%
\special{fp}%
\special{sh 1}%
\special{pa 3516 1450}%
\special{pa 3549 1511}%
\special{pa 3554 1488}%
\special{pa 3577 1483}%
\special{pa 3516 1450}%
\special{fp}%
\put(37.8000,-17.5000){\makebox(0,0)[lb]{$\psi$}}%
%
\special{pn 8}%
\special{pa 3400 680}%
\special{pa 3400 1610}%
\special{dt 0.045}%
\special{pa 3400 1610}%
\special{pa 3400 1609}%
\special{dt 0.045}%
\end{picture}%
\label{fig:D46}
\caption{Dynkin diagram of type $D_5^{(2)}$}
\end{figure}
\begin{equation}\label{19}
  \left\{
  \begin{aligned}
   \frac{dx}{dt} &=\frac{\partial H_{D_5^{(2)}}}{\partial y}=\frac{2x^2y-tx^2-2\alpha_0x+1}{t}-\frac{2x^2z(zw+\alpha_3)}{t},\\
   \frac{dy}{dt} &=-\frac{\partial H_{D_5^{(2)}}}{\partial x}=\frac{-2xy^2+2txy+2\alpha_0y+\alpha_1t}{t}+\frac{2z(zw+\alpha_3)(2xy+\alpha_1)}{t},\\
   \frac{dz}{dt} &=\frac{\partial H_{D_5^{(2)}}}{\partial w}=\frac{2z^2w-z^2+(1-2\alpha_4)z+t}{t}-\frac{2xz^2(xy+\alpha_1)}{t},\\
   \frac{dw}{dt} &=-\frac{\partial H_{D_5^{(2)}}}{\partial z}=\frac{-2zw^2+2zw-(1-2\alpha_4)w+\alpha_3}{t}+\frac{2x(xy+\alpha_1)(2zw+\alpha_3)}{t}
   \end{aligned}
  \right. 
\end{equation}
with the Hamiltonian \eqref{20}. Here $x,y,z$ and $w$ denote unknown complex variables and $\alpha_0,\alpha_1,ldots,\alpha_4$ are complex parameters satisfying the relation $\alpha_0+\alpha_1+\alpha_2+\alpha_3+\alpha_4=\frac{1}{2}$.

\begin{theorem}\label{4.1}
The system \eqref{19} admits extended affine Weyl group symmetry of type $D_5^{(2)}$ as the group of its B{\"a}cklund transformations {\rm (cf. \cite{3}), \rm} whose generators are explicitly given as follows{\rm : \rm}with the notation $(*):=(x,y,z,w,t;\alpha_0,\alpha_1,\ldots,\alpha_4)$,
\begin{align*}
        s_{0}: (*) &\rightarrow (-x,-y+\frac{2\alpha_0}{x}-\frac{1}{x^2},-z,-w,-t;-\alpha_0,\alpha_1+2\alpha_0,\alpha_2,\alpha_3,\alpha_4),\\
        s_{1}: (*) &\rightarrow (x+\frac{\alpha_1}{y},y,z,w,t;\alpha_0+\alpha_1,-\alpha_1,\alpha_2+\alpha_1,\alpha_3,\alpha_4), \\
        s_{2}: (*) &\rightarrow  (x,y-\frac{\alpha_2z}{xz-1},z,w-\frac{\alpha_2x}{xz-1},t;\alpha_0,\alpha_1+\alpha_2,-\alpha_2,\alpha_3+\alpha_2,\alpha_4), \\
        s_{3}: (*) &\rightarrow (x,y,z+\frac{\alpha_3}{w},w,t;\alpha_0,\alpha_1,\alpha_2+\alpha_3,-\alpha_3,\alpha_4+\alpha_3), \\
        s_{4}: (*) &\rightarrow (x,y,z,w-\frac{2\alpha_4}{w}+\frac{t}{z^2},-t;\alpha_0,\alpha_1,\alpha_2,\alpha_3+2\alpha_4,-\alpha_4), \\
        \psi: (*) &\rightarrow (\frac{z}{t},tw,tx,\frac{y}{t},t;\alpha_4,\alpha_3,\alpha_2,\alpha_1,\alpha_0).
\end{align*}
\end{theorem}

\begin{theorem}\label{4.2}
Let us consider a polynomial Hamiltonian system with Hamiltonian $H \in {\Bbb C}(t)[x,y,z,w]$. We assume that

$(A1)$ $deg(H)=5$ with respect to $x,y,z,w$.

$(A2)$ This system becomes again a polynomial Hamiltonian system in each coordinate $r_i \ (i=0,1,3,4)${\rm : \rm}
\begin{align*}
&r_0:x_0=x, \ y_0=y-\frac{2\alpha_0}{x}+\frac{1}{x^2}, \ z_0=z, \ w_0=w, \\
&r_1:x_1=1/x, \ y_1=-(yx+\alpha_1)x, \ z_1=z, \ w_1=w, \\
&r_3:x_3=x, \ y_3=y, \ z_3=1/z, \ w_3=-z(wz+\alpha_3), \\
&r_4:x_4=x, \ y_4=y, \ z_4=z, \ w_4=w-\frac{2\alpha_4}{w}+\frac{t}{z^2}.
\end{align*}

$(A3)$ In addition to the assumption $(A2)$, the Hamiltonian system in the coordinate $r_1$ becomes again a polynomial Hamiltonian system in the coordinate $r_2${\rm : \rm}
\begin{equation*}
r_2:x_2=-((x_1-z_1)y_1-\alpha_2)y_1, \ y_2=1/y_1, \ z_2=z_1, \ w_2=w_1+y_1.
\end{equation*}
Then such a system coincides with the system \eqref{19}.
\end{theorem}

Theorems \ref{4.1} and \ref{4.2} can be checked by a direct calculation, respectively.

\begin{theorem}\label{4.3}
For the system \eqref{1} of type $D_4^{(1)}$, we make the change of parameters and variables
\begin{gather}
\begin{gathered}\label{22}
A_0=\frac{\alpha_0-\alpha_1}{2}, \quad A_1=\alpha_1, \quad A_2=\alpha_2, \quad A_3=\alpha_3, \quad A_4=\frac{\alpha_4-\alpha_3}{2},\\
\end{gathered}\\
\begin{gathered}\label{23}
X=\frac{1}{x}, \quad Y=-(xy+\alpha_1)x, \quad Z=\frac{1}{z}, \quad W=-(zw+\alpha_3)z
\end{gathered}
\end{gather}
from $\alpha_0,\alpha_1, \dots ,\alpha_4,x,y,z,w$ to $A_0,A_1,\dots ,A_4,X,Y,Z,W$. Then the system \eqref{1} can also be written in the new variables $X,Y,Z,W$ and parameters $A_0,A_1,\dots ,A_4$ as a Hamiltonian system. This new system tends to the system \eqref{19} with the Hamiltonian \eqref{20}.
\end{theorem}

\begin{proof}
Notice that
$$
2(A_0+A_1+A_2+A_3+A_4)=\alpha_0+\alpha_1+2\alpha_2+\alpha_3+\alpha_4=1
$$
and the change of variables from $(x,y,z,w)$ to $(X,Y,Z,W)$ is symplectic. Choose $S_i \ (i=0,1,\dots ,4)$ and $\psi$ as
$$
S_0:=\pi_1, \ S_1:=s_1, \ S_2:=s_2, \ S_3:=s_3, \ S_4:=\pi_2, \ \psi:=\pi_3.
$$
Then the transformations $S_i$ are reflections of the parameters $A_0,A_1,\dots ,A_4$. The transformation group $\tilde{W}(D_5^{(2)})=<S_0,S_1,\dots ,S_4,\psi>$ coincides with the transformations given in Theorem \ref{4.1}.
\end{proof}

{\it Acknowledgement.} The author would like to thank K. Fuji, W. Rossman,  K. Takano and Y. Yamada for useful discussions.

\end{document}